# Naturalized bracket row and Motzkin triangle

Gennady Eremin
ergenns@gmail.com

April 20, 2020

**Abstract.** The paper considers the sequence of the Motzkin words, which is constructed according to formal features of natural numbers. We investigate the decomposition of well-formed parentheses into the matched pairs of parentheses (analogous to Prime Numbers). The point triangle of bracket pairs is constructed and studied, as well as the numerical Motzkin triangle (similar to Pascal's triangle). To calculate the Motzkin triangle, the graphical method of stars and bars is used.

The naturalized sequence of the Motzkin words is introduced and described in [Er19a, Er19b, Ere20]. The Motzkin words are totally ordered according to the formal features of natural numbers. The article continues to analyze the sequence, considers the decomposition of the Motzkin words into pairs of parentheses, and examines the Motzkin triangle – an analog of Pascal's triangle.

## 1  Introduction

This paper deals with well-formed parentheses strings with zeros, the *Motzkin words*, that are enumerated by Motzkin numbers [Wei19]. Interest in such parentheses is quite high. For example, in [BP14], the distance between the Motzkin words in the Tamari lattice is established. We investigate the total order in a naturalized sequence.

**Motzkin words**. There are different ways to represent the Motzkin words [DS77]. We work with bracket sets with sparse zeros in which (i) the number of left and right parentheses is the same and (ii) in any initial fragment of a word the number of left parentheses is no less than the right ones. A word can only contain zeros. Sets of length $n$, $n$-words, are counted by Motzkin number $M_n$, $n \geq 0$. The first Motzkin numbers are (see A001006): 1, 1, 2, 4, 9, 21, 51, 127, 323, 835, 2188, ... Here is the recurrent relation:

(1.1) $\qquad M_0 = 1, \ M_n = M_{n-1} + \sum_{k=0}^{n-2} M_k M_{n-2-k}, \ n \geq 1.$

The Motzkin number generating function $M(x)$, satisfies [Wei19]

(1.2) $\qquad M(x) = 1 + x^2 M^2(x) + xM(x) = \frac{1-x-\sqrt{1-2x-3x^2}}{2x^2} = \frac{2}{1-x+\sqrt{1-2x-3x^2}}.$

**Motzkin increments (difference numbers)**. The Motzkin *n*-word, $n \geq 1$, begins with either a *leading zero* or a left parenthesis. Words with a leading zero are inherited from the (*n*-1)-words by adding zero; the number of the *inherited words* is $M_{n-1}$. Other *n*-words, *unique words*, begin with a left parenthesis. The 1-word "0" is considered unique, and this is the only unique word with a leading zero. Based on (1.1), the number of the unique Motzkin *n*-words is equal to

(1.3) $\quad U_0 = 0, \ U_1 = 1, \ U_n = M_n - M_{n-1} = \sum_{k=0}^{n-2} M_k M_{n-2-k}, \ n \geq 2.$

The *increments* (1.3) will be called the *difference Motzkin numbers*. Obviously,

$$M_0 = 1, \ M_n = U_n + M_{n-1} = U_n + U_{n-1} + M_{n-2} = \ldots = \sum_0^n U_k, \ n \geq 1.$$

The Motzkin increments $U_n$ for $n = 0, 1, 2, \ldots$ form the sequence

(1.4) $\quad$ 0, 1, 1, 2, 5, 12, 30, 76, 196, 512, 1353, 3610, 9713, 26324, 71799, ….

From (1.3) one can calculate the generating function $U(x)$ [Ere20]:

(1.5) $\quad U(x) = x - 1 + (1-x)M(x). = x - 1 + \frac{2-2x}{1-x+\sqrt{1-2x-3x^2}}.$

**Natural numbers**. Many mathematicians consider zero to be a natural number, and among the Motzkin words there is also "0". Let $\mathbb{N}$ be a set of natural numbers and let $0 \in \mathbb{N}$. Natural numbers are indexed from 0, the index is equal to the value of the item, i.e. $n_i = i, i \geq 0$. Indexes are not repeated, and numbers are ordered by code length, by *ranges*. First there are single-digit items, 1-range, followed by two-digit numbers, 2-range, and so on. This is the *external sorting*.

In each range, numbers are sorted according to a strict order of alphabetic characters: $0 < 1 < \ldots < 9$. And this is the *internal sorting*. And last, an integer can begin with any alphabetic character, but not zero. The exception is the item zero itself.

*Note.* The set of natural numbers is self-indexing, the item is equal to its index, and a rhetorical question arises: when we work with an integer, maybe we are operating with its index? It is possible that we add, subtract, multiply (etc.) exactly indices. This does not matter for numbers, but it does for non-numeric sequences. In this connection, functions that operate on indexes in a non-numeric sequence are significant numerical functions that we can relate to functions that are defined on natural numbers.

**Naturalized Motzkin words.** Based on the formal properties of natural numbers, we get the totally ordered set of the unique Motzkin words, the *Motzkin row* (see [Er19b]):

$$\mathfrak{M} = \{0, \ (\ ), \ (0), \ (\ )0, \ (00), \ (0)0, \ ((\ )), \ (\ )00, \ (\ )(\ ), \ (000), \ (00)0, \ (0(\ )), \ldots\}$$

Items are indexed from zero and grouped by code length, *external sorting* forms *ranges*. In the ranges, we have *internal sorting*: the items are sorted according to the order in the alphabet 0 < ( < ). As in natural numbers, zero is a free sign, has a minimal weight, and is placed without restrictions. Parentheses are tied. In the Motzkin word (also in the Dick word), the brackets are divided into *matched pairs*. In $\mathfrak{M}$ and $\mathbb{N}$, the initial elements coincide and are equal to zero.

Obviously, $\#\mathfrak{M}_1 = U_1 = 1$, $\#\mathfrak{M}_n = U_n = M_n - M_{n-1}$, $n > 1$. We work with real bracket sets, so $\#\mathfrak{M}_0 = U_0 = 0$. And this is logical; there is no empty item in $\mathbb{N}$. Each range has a minimum and a maximum:

$$\min \mathfrak{M}_n = (0^{n-2}) \quad \text{and} \quad \max \mathfrak{M}_n = ()^k[0], \text{ where } k = \lfloor n/2 \rfloor.$$

In this case, $\lfloor x \rfloor$ is the floor function, the greatest integer less than or equal to $x$. The superscript indicates the repetition of zero and adjacent parentheses. The maximum ends with 0 if $n$ is odd. Let's call the index of $x \in \mathfrak{M}$ *weight* and denote wt $x$. So

$$\text{wt }(0)0 = 5, \quad \text{wt}(\min \mathfrak{M}_n) = M_{n-1}, \quad \text{wt}(\max \mathfrak{M}_n) = \text{wt}(\min \mathfrak{M}_{n+1}) - 1 = M_n - 1.$$

**Block arithmetic.** A matched pair and everything inside is called a *block*. There is an *internal block* and an *external block*, aka a *simple block* or a *simple word*. The internal block is nested in another block. A simple word is not contained in any block, i.e. it is at the *zero level* of bracket nesting.

A simple blocks can be added (superimposed) or extracted from bracket sets, and the positions of the brackets should be fixed [Er19a]. A simple block may acquire additional trailing zeros during the extraction process. For example, from $()0(())(00)$ you can extract the internal block $(())0^4$, write it like this: $()0(())(00) = ()0^5(00) + (())0^4$, or using the weights: wt $()0(())(00) = $ wt $()0^5(00) + $ wt $(())0^4$.

**A little *p*-adica.** We can consider the Motzkin word as a set of brackets that are superimposed on an *infinite zero line* bounded on the right. The Motzkin word adjoins the border of such a line. You can move the brackets to the left unlimitedly, increasing the range of the word and adding trailing zeros. But moving to the right is impossible if there are no trailing zeros in the word. In the case of processing a group of the Motzkin words such an infinite line can serve as an *adder*.

We see a similar interpretation in the *p*-adic analysis when we work with binary integers and when an infinite string of 0's and 1's is limited to the right [Ana14].

## 2    Matched pairs of parentheses

Let's call a block with one pair of brackets a *prime pair*. In the general case, the Motzkin word consists of a set of matched pairs that are spaced at different levels of nesting, starting from level zero, the level of simple blocks. The Motzkin word is decomposed into prime pairs (how not to recall the decomposition of an integer into a product of primes?). Here is a set of prime pairs:

$$\mathscr{P} = \{\,(\,),\,(0),\,(\,)0,\,(00),\,(0)0,\,(\,)00,\,(0^3),\,(00)0,\,(0)00,\,(\,)0^3,\,(0^4),\,\ldots\} \subset \mathfrak{M}.$$

It is convenient to index prime pairs starting from 1: $p_1 = (\,)$, $p_2 = (0)$, and so on, since the initial elements get the same indexes in both $\mathcal{P}$ and $\mathfrak{M}$. The elements are distributed in ranges. In the *n*-range there is the smallest pair $(0^{n-2})$ without trailing zeros and the largest pair $(\,)0^{n-2}$ with zeros in the tail. Obviously, $\#\mathcal{P}_n = n-1$.

We denote by $p_{n,r}$ the prime pair of length $n$ and with a right parenthesis in the position $r$, so $p_{n,r} = (0^{n-r-1})0^{r-1}$, $n > r \geq 1$. Obviously, $p_{n,1} = (0^{n-2})$, $p_{n,n-1} = (\,)0^{n-2}$. It is easy to see, $p_{n,r} = p_i \in \mathcal{P}$, $i = r + T_{n-2}$, where $T_{n-2} = 1 + 2 + \ldots + (n-2) = (n-1)(n-2)/2$ (see the *triangular numbers* [Wik20b, A000217]).

**Weight function of a prime pair.** Let's calculate the weight of a prime pair in $\mathfrak{M}$. Table 1 shows the first 60 items from $\mathfrak{M}$. The left column is supplemented with weights; the cells with prime pairs are darkened. The last marked pair $p_{7,4} = (00)0^3$ has a weight of 59.

| | | | | | | | | | |
|---|---|---|---|---|---|---|---|---|---|
| 0: 0 | (\,) | (0) | (\,)0 | (00) | (0)0 | ((\,)) | (\,)00 | (\,)(\,) | $(0^3)$ |
| 10: (00)0 | (0(\,)) | (0)00 | (0)(\,) | ((0)) | ((\,)0) | ((\,))0 | (\,)$0^3$ | (\,)0(\,) | (\,)(0) |
| 20: (\,)(\,)0 | $(0^4)$ | $(0^3)0$ | (00(\,)) | (00)00 | (00)(\,) | (0(0)) | (0(\,)0) | (0(\,))0 | (0)$0^3$ |
| 30: (0)0(\,) | (0)(0) | (0)(\,)0 | ((00)) | ((0)0) | ((0))0 | (((\,))) | ((\,)00) | ((\,)0)0 | ((\,)(\,)) |
| 40: ((\,))00 | ((\,))(\,) | (\,)$0^4$ | (\,)000(\,) | (\,)0(0) | (\,)0(\,)0 | (\,)(00) | (\,)(0)0 | (\,)(\,)(\,) | (\,)(\,)00 |
| 50: (\,)(\,)0(\,) | $(0^5)$ | $(0^4)0$ | $(0^3)(\,)$ | $(0^3)00$ | $(0^3)(\,)$ | (00(0)) | (00(\,)0) | (00(\,))0 | (00)$0^3$ |

Table 1: First prime pairs (dark cells).

Obviously, wt $p_{n,1} = M_{n-1}$, wt $p_{n+1,1} = M_n$. Moving the left parenthesis from position $n$ to position $n+1$ (moving a prime pair from range $n$ to range $n+1$) increases the weight by $M_n - M_{n-1}$. The increment is independent of the position of the right parenthesis. Let's find a general formula for wt $p_{n,r}$.

We begin the calculations with the pair $p_{r+1,r} = (\,)0^{r-1}$, the weight of this pair is

$$\text{wt } p_{r+1,r} = \text{wt}(\max \mathfrak{M}_{r+1}) - \text{wt}(\max \mathfrak{M}_{r-1}) = M_{r+1} - M_{r-1}.$$

Let's increase the weight by moving the left parenthesis to position $n$ step by step:

$p_{r+1,r} = (\,)0^{r-1}$ : wt $p_{r+1,r} = \phantom{M_{r+1} - M_r +{}} M_{r+1} - M_{r-1}$;
$p_{r+2,r} = (0)0^{r-1}$ : wt $p_{r+2,r} = M_{r+1} - M_r + M_{r+1} - M_{r-1}$;
$p_{r+3,r} = (00)0^{r-1}$ : wt $p_{r+3,r} = M_{r+2} - M_r + M_{r+1} - M_{r-1}$;
. . . . . . . . .
$p_{n,r} = (0^{n-r-1})0^{r-1}$ : wt $p_{n,r} = M_{n-1} - M_r + M_{r+1} - M_{r-1}$.

As a result, we get

(2.1) $\qquad\qquad \text{wt } p_{n,r} = M_{n-1} + M_{r+1} - M_r - M_{r-1}, \ n > r \geq 1.$

Check the weight of the pair $p_{7,4}$: wt $p_{7,4} = M_6 + M_5 - M_4 - M_3 = 51 + 21 - 9 - 4 = 59$. Next we separate from (2.1) the *r*-part:

(2.2) $\qquad\qquad \text{wt } p_{n,r} = M_{n-1} + \delta_r, \ \delta_r = U_{r+1} - M_{r-1}, \ n > r \geq 1.$

The numbers $\delta_r$ of (2.2) form the sequence for $r = 1, 2, \ldots$ :

(2.3)     0, 1, 3, 8, 21, 55, 145, 385, 1030, 2775, 7525, 20526, 56288 …

**Generating function of a prime pair.** The weight function (2.1) depends on two parameters, i.e. we have a two-index sequence. Such sequences are often written as a triangle (for example, Pascal's triangle). Table 2 shows the weights of prime pairs. The first column shows the length of prime pairs (the range or the position of the left parenthesis) starting with $n = 2$; on the top line are the positions of the right parenthesis starting with $r = 1$. The $n$th line contains the weights of prime pairs of the $n$-range, the weight increments correspond to the sequence (2.3). We highlighted in red the weight of the aforementioned pair $p_{7,4}$.

| n\r | 1 | 2 | 3 | 4 | 5 | 6 | 7 |
|---|---|---|---|---|---|---|---|
| 2 | 1 | | | | | | |
| 3 | 2 | 3 | | | | | |
| 4 | 4 | 5 | 7 | | | | |
| 5 | 9 | 10 | 12 | 17 | | | |
| 6 | 21 | 22 | 24 | 29 | 42 | | |
| 7 | 51 | 52 | 54 | 59 | 72 | 106 | |
| 8 | 127 | 128 | 130 | 135 | 148 | 182 | 272 |

Table 2: Triangle of prime pair weights.

According to (2.2), we write the generating function as follows:

(3.4)     $\text{Pair}(x, y) = \sum_{n \geq 2} M_{n-1} x^n + \sum_{r \geq 1} (U_{r+1} - M_{r-1}) y^r$.

First we calculate the $n$-part of the generating function:

$\sum_{n \geq 2} M_{n-1} x^n = x \sum_{n \geq 2} M_{n-1} x^{n-1} = x(M(x) - 1)$.

Next, we find the $r$-part:

$\sum_{r \geq 1} (U_{r+1} - M_{r-1}) y^r = \sum_{r \geq 1} U_{r+1} y^r - \sum_{r \geq 1} M_{r-1} y^r$

$= \frac{1}{y} \sum_{r \geq 1} U_{r+1} y^{r+1} - y \sum_{r \geq 1} M_{r-1} y^{r-1}$

$= \frac{1}{y} (U_2 y^2 + U_3 y^3 + \ldots + U_0 + U_1 y - U_0 - U_1 y) - yM(y)$

$= \frac{1}{y} (U(y) - y) - yM(y) = \frac{1}{y} U(y) - 1 - yM(y)$

$= \frac{1}{y} (y - 1 + (1 - y)M(y)) - 1 - yM(y) = \frac{1-y-y^2}{y} M(y) - \frac{1}{y}$.

We used (1.5), as a result we get

$\text{Pair}(x, y) = xM(x) - x + \frac{1-y-y^2}{y} M(y) - \frac{1}{y}$

$$= \frac{1-x-\sqrt{1-2x-3x^2}}{2x} - x + \frac{2-2y-2y^2}{y(1-y+\sqrt{1-2y-3y^2})} - \frac{1}{y}.$$

## 3  Motzkin triangle

**Dot triangle**. Let's rotate Table 2 counterclockwise and represent the prime pairs on the plane as dots (see Figure 1). On the abscissa axis, we postponed the length of prime pairs, the positions of the right parentheses are indicated on the ordinate axis. All points are located between the abscissa axis and the diagonal $r = n$. Each prime pair is represented by a point in such a triangle and vice versa. Some prime pairs (including $p_{7,4}$) are represented by weights from Table 2.

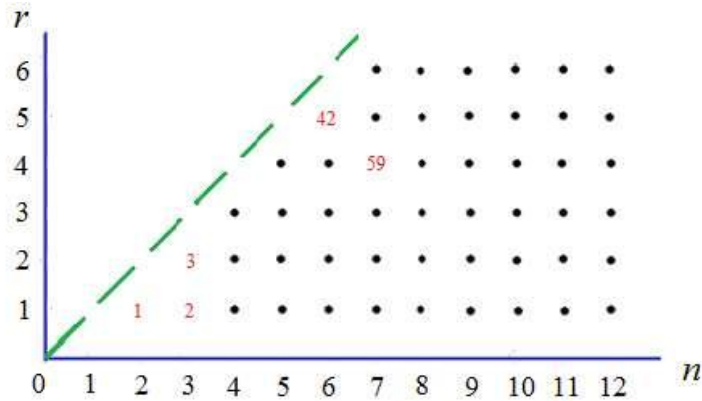

Figure 1: Dot triangle of prime pairs.

It is easy to see that the $n$-range from $\mathfrak{M}$ is bounded by the column with $n-1$ points: the lower point $(n, 1)$ is the minimal prime pair $p_{n,1}$, and the upper point $(n, n-1)$ is the maximal pair $p_{n,n-1}$ in the $n$-range. Any prime pair from such an *n-triangle* (bounded by points with abscissa $n$) can appear inside the Motzkin $n$-word.

Well-parenthesed word is uniquely decomposed into prime pairs, so the Motzkin word corresponds to a certain set of points in Figure 1. For example, three points with weights 3, 42, and 59 correspond to the decomposition of the word (())()0 into prime pairs $p_{3,2} = ()0$, $p_{6,5} = ()0000$, and $p_{7,4} = (00)000$. We can calculate the weight of (())()0  (see [Er19b]):

$$\text{wt }(())()0 \ = \ \text{wt } p_{3,2} + \text{wt}' p_{6,5} + \ \text{wt } p_{7,4} \ = 3 + 43 + 59 = 105.$$

In the Motzkin word, the leftmost parenthesis belongs to a prime pair that defines the range of the word. This pair can be called the *main* one. It is more convenient to start the decomposition of the Motzkin word with the main pair.

Let $z$ be the Motzkin word and let's choose the prime pair $p_{a,b}$, $a > b \geq 1$, from $z$. Figure 2 shows the pair $p_{a,b}$ as a point with coordinates $a, b$. We have identified six areas (zones) of the dot triangle: to the left of the point $(a, b)$ there are three fi-nite areas bounded by the diagonal and the abscissa axis; on the right there are three infinite areas. Let's explain the coloring of the zones.

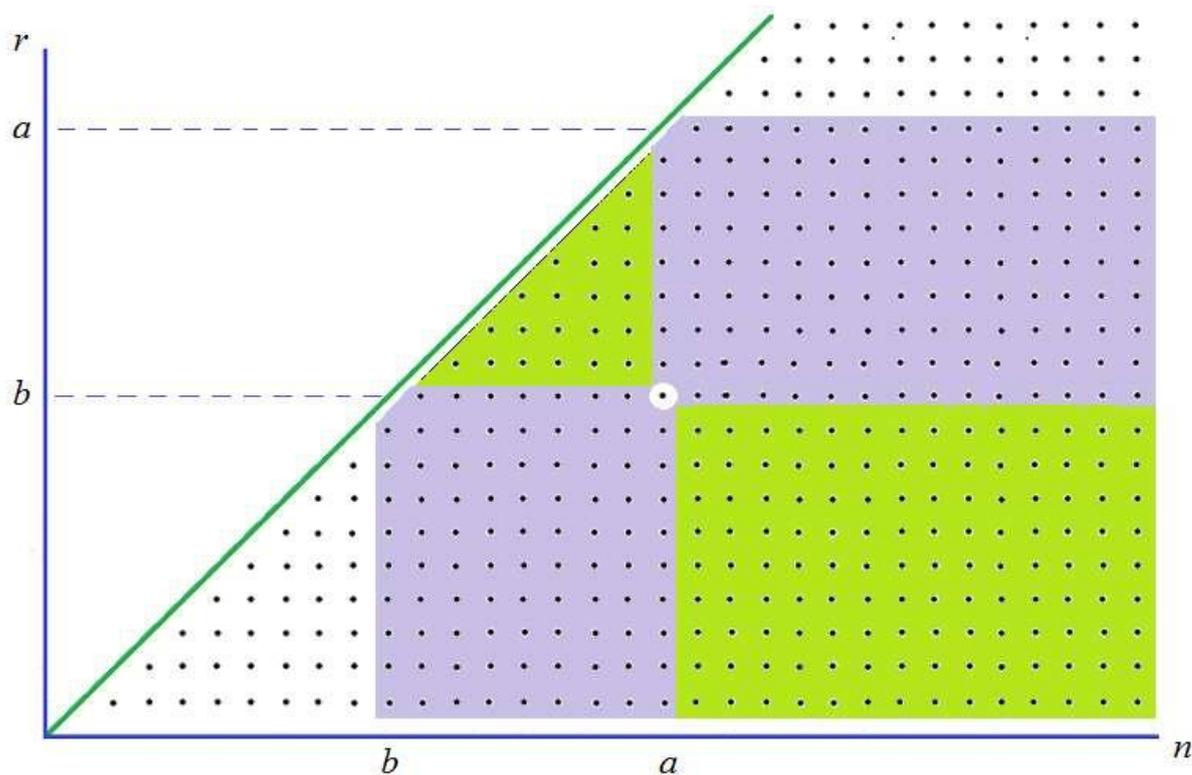

Figure 2: Zones of the dot triangle.

Two zones are darkened; these points are incompatible with the prime pair $p_{a,b}$. The parentheses of each pair occupy fixed positions; we can say that the point $(a, b)$ has "reserved" the positions $a, b$ for itself. Neither the left nor right parentheses of other pairs in $z$ can occupy these positions. In addition, parentheses cannot intersect, for example, the adjacent point $(a-1, b-1)$ from the smaller area (lower gray zone) and the nearest pair $(a+1, b+1)$ from the larger area (upper gray zone) intersect with the given point $(a, b)$.

In other zones, the points are compatible with the pair $(a, b)$. In green, we have shown two areas where the points are either nested inside $p_{a,b}$ (a bounded zone near the diagonal), or they themselves include this pair (an infinite zone along the abscissa axis). The two zones are not painted, and the points of these zones do not intersect with $p_{a,b}$. We are interested in the zone at the origin where the initial elements of $\mathscr{P}$ are located.

**$n$-triangle.** Further we will work in the $n$-range of $\mathfrak{M}$. Figure 3 shows three zones for the point $(n, r)$, $n > r \geq 1$. We shaded the closed area again; these points cannot appear in the same Motzkin word together with $(n, r)$. While the points of the triangles A and B are compatible with $(n, r)$: the points of A do not intersect with the given pair, and the points of B are nested in this pair.

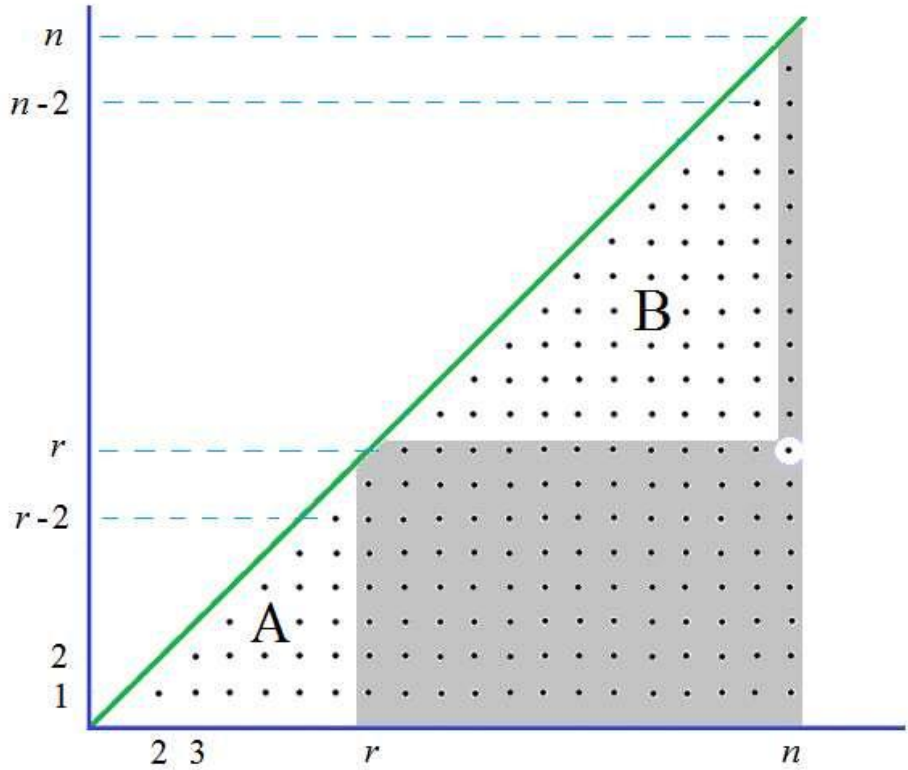

Figure 3: Dot $n$-triangle.

We denote by $\mathfrak{M}_{n,k}$ set of the Motzkin $n$-word that contains exactly $k$ prime pairs. Accordingly, $U_{n,k} = \#\mathfrak{M}_{n,k}$. Obviously, for $n \geq 2$ $\mathfrak{M}_n = \bigcup_k \mathfrak{M}_{n,k}$, $U_n = \sum_k U_{n,k}$.

The number of the Motzkin $n$-words with one pair of parentheses is

(3.1) $$U_{n,1} = n-1, \ n \geq 1,$$

and exactly so many points are in the $n$-column of the dot triangle.

Let's calculate the number of $n$-words with two pairs of parentheses, $U_{n,2}$. We need to count the number of points in triangles A and B if, in the pair $(n, r)$, we iterate over all the positions of the right parenthesis $r = 1, 2, \ldots n-1$. In A and B, the dot configuration is similar to a regular triangle (the number of points on the sides is the same). Let's look at a few cases.

*Case $r = 1$.* The closed zona is the bottom line and the $n$th column, so zone A is empty, $\#A = 0$. But zone B has the maximum number of points, $\#B = 1 + 2 + \ldots + n-3 = T_{n-3}$ (a triangular number appears again [Wi20b]).
*Symmetrically $r = n-1$.* The closed zone is columns $n-1$ and $n$. Now the picture is reversed: $\#A = T_{n-3}$, and zone B is degenerating, $\#B = 0$.

*Case $r = 2$.* The closed zone occupies the bottom two lines and the $n$th column. Zone A is still empty, $\#A = 0$, and B has one less line, $\#B = T_{n-4}$.
*Symmetrically $r = n-2$.* In the closed zone we have the last three columns. The image is inverted again: in A, there is one line less, $\#A = T_{n-4}$, and $\#B = 0$.

*Case r* = 3. The first point (2, 1) appears in zone A, i.e. #A = 1 = $T_1$, in B again there is one less line, #B = $T_{n-5}$.

*Symmetrically r* = *n*–3. In the closed zone, we have the last three columns without a top point. Again the picture is inverted: #A = $T_{n-5}$, and in B we a single point (*n*–1, *n*–2), т.е. #B = $T_1$.

Thus, when moving the pair (*n*, *r*) along the *n*th column, the zones A and B take the same dimensions, but in the reverse order. Therefore, we only need to calculate one zone and double the result. It is more convenient to work with zone A, for example, in Figure 3, we simply count the points from the origin: #A = 1 + 2 + … + *r*–2 = $T_{r-2}$, *r* > 2.

The considered cases indicate the need to sum the triangular numbers from 1 to $T_{n-3}$ (then we double the resulting value). In mathematics, the sum of the first *t* triangular numbers denotes $Te_t$ and is called the *t*-th *tetrahedral number* [Wei20]. The following formula is known:

$$Te_t = T_1 + T_2 + \ldots + T_t = \binom{t+2}{3} = t(t+1)(t+2)/6.$$

**Proposition 3.1.** *The number of the Motzkin n-word with two pairs of parentheses is*

(3.2) $\qquad U_{n,2} = 2Te_{n-3} = 2\binom{n-1}{3} = (n-1)(n-2)(n-3)/3, \quad n > 3.$

The numbers $U_{n,2}$ for *n* = 4, 5, 6, … are doubled elements of the sequence A000292: 2, 8, 20, 40, 70, 112, 168, 240, 330, 440, ... In Table 3, we highlighted prime pairs in red, and in green, the Motzkin words with two pairs of parentheses. As you can see, in the 6-range there are 5 prime pairs and 20 items with two pairs of parentheses. But there are still 5 words with three pairs of parentheses.

| 0: | 0 | () | (0) | ()0 | (00) | (0)0 | (()) | ()00 | ()() | $(0^3)$ |
|---|---|---|---|---|---|---|---|---|---|---|
| 10: | (00)0 | (0()) | (0)00 | (0)() | ((0)) | (()0) | (())0 | $()0^3$ | ()0() | ()(0) |
| 20: | ()()0 | $(0^4)$ | $(0^3)0$ | (00()) | (00)00 | (00)() | (0(0)) | (0()0) | (0())0 | $(0)0^3$ |
| 30: | (0)0() | (0)(0) | (0)()0 | ((00)) | ((0)0) | ((0))0 | ((())) | (()00) | (()0)0 | (()() ) |
| 40: | (())00 | (())() | $()0^4$ | ()00() | ()0(0) | ()0()0 | ()(00) | ()(0)0 | ()(()) | ()()00 |
| 50: | ()()() | $(0^5)$ | $(0^4)0$ | $(0^3)()$ | $(0^3)00$ | $(0^3)()$ | (00(0)) | (00()0) | (00())0 | $(00)0^3$ |

Table 3: Distribution of the Motzkin words with 1 and 2 pairs of parentheses.

**Numerical Motzkin triangle.** In Table 4 below, we have shown the results of direct calculation of $U_{n,k}$ for the first 15 ranges with the number of prime pairs up to 7 inclusive, zero elements are removed. The lines indicate the number of *n*-words with a fixed number of prime pairs, and the first two lines are already familiar to us.

The sum of the numbers in the *n*th column is $U_n$. It is logical to call Table 4 the *Motzkin triangle*. Let's formulate a statement that lists some properties of the Motzkin triangle.

| k\n | 2 | 3 | 4 | 5 | 6 | 7 | 8 | 9 | 10 | 11 | 12 | 13 | 14 | 15 |
|---|---|---|---|---|---|---|---|---|---|---|---|---|---|---|
| 1 | 1 | 2 | 3 | 4 | 5 | 6 | 7 | 8 | 9 | 10 | 11 | 12 | 13 | 14 |
| 2 | | | 2 | 8 | 20 | 40 | 70 | 112 | 168 | 240 | 330 | 440 | 572 | 728 |
| 3 | | | | | 5 | 30 | 105 | 280 | 630 | 1260 | 2310 | 3960 | 6435 | 10010 |
| 4 | | | | | | | 14 | 112 | 504 | 1680 | 4620 | 11088 | 24024 | 48048 |
| 5 | | | | | | | | | 42 | 420 | 2310 | 9240 | 30030 | 84084 |
| 6 | | | | | | | | | | | 132 | 1584 | 10296 | 48048 |
| 7 | | | | | | | | | | | | | 429 | 6006 |
| ... | | | | | | | | | | | | | | |
| ∑ | 1 | 2 | 5 | 12 | 30 | 76 | 196 | 512 | 1353 | 3610 | 9713 | 26324 | 71799 | 196938 |

Table 4: Motzkin triangle.

**Proposition 3.2**. *For the k-th line of the Motzkin triangle*
  (i)  $U_{<2k,k} = 0$,
  (ii) $U_{2k,k} = Cat(k)$,
  (iii) $U_{2k+1,k} = 2k \times Cat(k)$,
  (iv) $U_{2k+2,k} = k(2k+1) \times Cat(k)$.

*Proof.* Let's look at each item separately.

  (i) The Motzkin *n*-words, $n < 2k$, cannot contain $k$ or more prime pairs. This is why Table 4 is triangular.

  (ii) The Motzkin (2k)-words can contain $k$ prime pairs, but in this case we are dealing with a Dyck word (without zeros); the number of such sets is equal to the *k*th Catalan number, $Cat(k)$.

  (iii) Let's add a zero to the Dyck of length 2k; in this case, we obtain the Motzkin (2k+1)-word. We can add zero after any of the 2k parentheses, so the number of such Motzkin words is $2k \times Cat(k)$.

  (iv) If two zeros are added to the Dyck word of length 2k, then we get the Motzkin word from the (2k+2)-range. How many ways can we add two zeros? For placement, we have 2k places (after each parenthesis), and you can insert both zeros in each place at once. As a result, we get 2k placement options (as in the previous variant). But zeros can be placed one at a time in different places; if one zero is placed in one of the 2k places, then there are 2k–1 places for the second zero. In total, we get another 2k(2k–1) ways, but this number must be divided by 2 to remove duplicates (zeros are identical). Thus, the number of ways is equal to
$U_{2k+2,k} = Cat(k) \times (2k + 2k(2k-1)/2) = k(2k+1)\,Cat(k)$  □

Let's find a general formula for calculating $U_{n,k}$, $n \geq 2k$, use graphical method that William Feller suggested [Fel68, Wi20a].

**Stars and bars method.** This method is used for solving some problems, such as how many ways there are to put indistinguishable balls (stars) into distinguishable boxes. During the analysis, the balls and partitions (bars) between the boxes are lined up; balls can be placed in front of the first partition, between partitions, and behind the last partition. In our case, the boxes can be empty, i.e. there may be several partitions between adjacent balls.

In the Motzkin word, zeros are balls, and after each parenthesis we have a conditional box. If you delete the first parenthesis in the Motzkin word, then the remaining brackets can be considered partitions. Let's look at an example.

**Example** 3.3. In the word ((00())0)0, there are three pairs, $k = 3$, and four zeros. The word length is $n = 2k + 4 = 10$. Let's calculate $U_{10, 3}$. If we remove the first bracket and replace the remaining ones with partitions, we get the tuple |oo|||o|o of length $n - 1 = 9$ and with the number of partitions $2k - 1 = 5$. Feller's method is reduced to determining the number of permutations of partitions in the resulting tuple, in other words, we just need to calculate the following value: $\binom{n-1}{2k-1} = \binom{9}{5} = 9!/(5!4!) = 126$. It remains to multiply the obtained result by the Catalan number $Cat(3) = 5$.

As a result, we get $U_{10, 3} = \binom{10-1}{2*3-1} Cat(3) = 630$. Obviously, with the same result, it was possible to count the number of permutations of zeros in the tuple |oo|||o|o. □

Now we can state the corresponding lemma.

**Lemma 3.4**. *The number of the Motzkin n-words with k pairs of parentheses is*

(3.3) $$U_{n,k} = \binom{n-1}{2k-1} Cat(k) = \binom{n-1}{2k-1} \binom{2k}{k}/(k+1), \quad n \geq 2k > 1.$$

Below, we inverted Table 4 around the diagonal, giving it the usual form of Pascal's triangle (for example, see the Motzkin triangle in [DS77], p. 293).

| n\k | 1 | 2 | 3 | 4 | 5 | 6 | 7 | 8 | 9 |
|---|---|---|---|---|---|---|---|---|---|
| 2 | 1 | | | | | | | | |
| 3 | 2 | | | | | | | | |
| 4 | 3 | 2 | | | | | | | |
| 5 | 4 | 8 | | | | | | | |
| 6 | 5 | 20 | 5 | | | | | | |
| 7 | 6 | 40 | 30 | | | | | | |
| 8 | 7 | 70 | 105 | 14 | | | | | |
| 9 | 8 | 112 | 280 | 112 | | | | | |
| 10 | 9 | 168 | 630 | 504 | 42 | | | | |
| 11 | 10 | 240 | 1260 | 168 | 420 | | | | |
| 12 | 11 | 330 | 2310 | 4620 | 2310 | 132 | | | |
| 13 | 12 | 440 | 3960 | 11088 | 9240 | 1584 | | | |
| 14 | 13 | 572 | 6435 | 24024 | 30030 | 10296 | 429 | | |
| 15 | 14 | 728 | 10010 | 48048 | 84084 | 48048 | 6006 | | |
| 16 | 15 | 910 | 15015 | 90090 | 210210 | 180180 | 45045 | 1430 | |
| 17 | 16 | 1120 | 21840 | 160160 | 480480 | 576576 | 240240 | 22880 | |
| 18 | 17 | 1360 | 30940 | 272272 | 1021020 | 1633632 | 1021020 | 194480 | 4862 |

Table 4a: Motzkin triangle.

In the *n*th line, the first element is $n - 1$. In an even line $n = 2k$, the last element is $Cat(k)$, and in the case of $n = 2k + 1$, the last element is $2k\, Cat(k)$. In Table 4a, we have shown in red *matching elements*, the duplicates in lines. In this regard, we will formulate another statement.

**Proposition 3.5**. *In table 4a, in lines n = 3k, k > 1, the elements in columns k–1 and k+1 are equal, i. e.* $U_{3k,\,k-1} = U_{3k,\,k+1}$.

*Proof.* Let's use the formula (3.3).

$$U_{3k,\,k-1} = \binom{3k-1}{2k-3}\binom{2k-2}{k-1}/k = \frac{(3k-1)!\,(2k-2)!}{(2k-3)!\,(k+2)!\,(k-1)!\,k!}$$

$$= \frac{(3k-1)!\,(2k-2)}{(k+2)!\,(k-1)!\,k!} = \frac{2(3k-1)!}{(k+2)!\,(k-2)!\,k!};$$

$$U_{3k,\,k+1} = \binom{3k-1}{2k+1}\binom{2k+2}{k+1}/(k+2) = \frac{(3k-1)!\,(2k+2)!}{(2k+1)!\,(k-2)!\,(k+1)!\,(k+2)!}$$

$$= \frac{(3k-1)!\,(2k+2)}{(k-2)!\,(k+1)!\,(k+2)!} = \frac{2(3k-1)!}{(k-2)!\,k!\,(k+2)!} = U_{3k,\,k-1}.$$

□

In the line $n = 3k$, the maximum is located in the column $k$. It is easy to get equality: $(k-1)\,U_{3k,\,k} = (k+2)\,U_{3k,\,k+1}$, $k > 1$. Then the sum of the three maximum numbers in the line is $U_{3k,\,k-1} + U_{3k,\,k} + U_{3k,\,k+1} = \frac{3k}{k+2}\,U_{3k,\,k}$.

**Remark**. Formula (3.3) makes it possible to specify the amount in (1.3):

(3.4) $$U_n = \sum_{k=1}^{\lfloor n/2 \rfloor} \binom{n-1}{2k-1} Cat(k),\ n \geq 2.$$

In (3.4), the summands acquire a specific meaning – the number of Motzkin words with a fixed number of bracket pairs.

# References

[Ana14]  V. S. Anashin. *The p-adic ergodic theory and applications,* 2014. https://www.researchgate.net/publication/269571423

[BP14]   J.-L. Baril and J.-M. Pallo. *Motzkin subposets and Motzkin geodesics in Tamari lattices.* Inform. Process. Lett. 114 (1-2), 31-37, 2014. http://jl.baril.u-bourgogne.fr/Motzkin.pdf

[DS77]   R. Donaghey and L. W. Shapiro. *Motzkin numbers.* Journal of Combinatorial Theory, Series A, 23 (3), 291–301, 1977. https://www.sciencedirect.com/science/article/pii/0097316577900206

[Er19a]  G. Eremin, *Arithmetic on Balanced Parentheses*: *The case of Ordered Motzkin Words*, 2019. arXiv:1911.01673

[Er19b]  G. Eremin. *Mathematics of Balanced Parentheses*: *The case of Ordered Motzkin Words*, 2019. arXiv:1912.09693

[Ere20]  G. Eremin. *Generating function for Naturalized Series*: *The case of Ordered Motzkin Words*, 2019. arXiv:2002.08067

[Fel68]   W. Feller. *An Introduction to Probability Theory and Its Applications,* John Wiley & Sons, Inc., New York, 1950 (first edition).

[FS09]    Philippe Flajolet and Robert Sedgewick. *Analytic Combinatorics*. http://algo.inria.fr/flajolet/Publications/book.pdf

[Slo19]   N. J. A. Sloane. *The on-line encyclopedia of integer sequences*, 2019. https://oeis.org/

[Wei19]   Eric W. Weisstein. *Motzkin Number*. A Wolfram Web Resource, 2019. http://mathworld.wolfram.com/MotzkinNumber.html

[Wei20]   Eric W. Weisstein. *Tetrahedral Number*. A Wolfram Web Resource,. http://mathworld.wolfram.com/TetrahedralNumber.html

[Wi20a]   Stars and bars (combinatorics). https://en.wikipedia.org/wiki/Stars_and_bars_(combinatorics)

[Wi20b]   Triangular number. https://en.wikipedia.org/wiki/Triangular_number

Gzhel State University, Moscow, 140155, Russia
http://www.en.art-gzhel.ru/